\documentclass[12pt]{article}
\usepackage{fullpage}

\title{New  Stability Conditions for  Linear
Differential Equations with Several Delays}
\author{Leonid Berezansky $^{1}$ \\Department of Mathematics , 
Ben-Gurion University of the Negev, \\ 
Beer-Sheva 84105, Israel
\\
and  Elena Braverman $^{2}$ \\
Department of Mathematics and Statistics, University of Calgary, \\ 
2500 University Drive N.W., Calgary, AB T2N 1N4, Canada}
\date{}
\begin{document}
\maketitle

\footnotetext[1]{Partially supported by Israeli Ministry of Absorption}
\footnotetext[2]{Partially supported by the NSERC Research Grant}

\newcommand{\qed}{\hbox to 0pt{}\hfill$\rlap{$\sqcap$}\sqcup$ 
\vspace{2mm}}

\begin{abstract}
New explicit conditions of asymptotic and exponential stability 
are obtained for the scalar nonautonomous linear delay differential 
equation 
$$
\dot{x}(t)+\sum_{k=1}^m  a_k(t)x(h_k(t))=0
$$
with measurable delays and coefficients. 
These results are compared to known stability tests.
\end{abstract}


{\bf Running title:} New Stability Conditions

\thispagestyle{empty}

\section{Introduction}
\label{introduct}

In this paper we  continue the study of stability  properties
for the linear differential equation with several delays and an  arbitrary
number of positive and negative coefficients
\begin{equation}
\label{1}
\dot{x}(t)+\sum_{k=1}^m a_k(t)x(h_k(t))=0, ~~t \geq t_0,
\end{equation}
which was begun in \cite{BB1}-\cite{BB3}.
Equation  (\ref{1}) and its special
cases were intensively studied, for example, in  
\cite{YS1}-\cite{Mal2}.
In \cite{BB2} we gave a review of stability tests obtained in these
papers. 

In almost all papers on stability of delay differential equations
coefficients and delays are assumed to be continuous, which is essentially
used in the proofs of main results.
In real world problems, for example, in biological and ecological models 
with seasonal fluctuations of parameters and in economical models with
investments, parameters of differential equations are not necessarily 
continuous.

There are also some mathematical reasons to consider differential equations
without the assumption that parameters are continuous functions.
One of the main methods to investigate impulsive differential equations 
is their reduction to a non-impulsive differential equation with 
discontinuous coefficients.
Similarly, difference equations can sometimes be reduced to the similar 
problems for delay differential equations with discontinuous
piecewise constant  delays.

In paper \cite{BB1}  some problems for differential equations 
with several delays were reduced to
similar problems for equations with one delay which 
generally is not continuous.

One of the purposes of this paper is to extend and partially  improve most popular
stability results for linear delay equations with continuous coefficients and delays
to equations with measurable parameters.

Another purpose is to generalize some results of \cite{BB1,BB2, 
BB3}. In these papers, the sum of coefficients was supposed to be 
separated from zero and delays were
assumed to be bounded. So the results of these papers
are not applicable, for
example, to the following   equations
$$
\dot{x}(t)+|\sin t| x(t-\tau)=0,
$$$$
\dot{x}(t)+(|\sin t|-\sin t) x(t-\tau)=0,
$$$$
\dot{x}(t)+\frac{1}{t}x(t)+\frac{\alpha}{t}
x \left( \frac{t}{2} \right) =0.
$$
In most results of the present paper
these restrictions are omitted, so we can consider all the equations
mentioned above.
Besides, necessary stability conditions (probably for the first time) 
are obtained for equation (\ref{1}) with nonnegative coefficients and 
bounded delays.
In particular, if this equation is exponentially stable then the ordinary 
differential equation
$$
\dot{x}(t)+\sum_{k=1}^m a_k(t)x(t)=0
$$
is also exponentially stable.

\section{Preliminaries}
\label{prelim}

We consider the scalar linear equation with several 
delays (\ref{1}) for $t\geq t_0$ with the initial conditions 
(for any $t_0\geq 0$) 
\begin{equation}
\label{2}
x(t)=\varphi(t), ~ t<t_0,~ x(t_0)=x_0, 
\end{equation}
and under the following assumptions:

(a1) $a_k(t)$ are Lebesgue measurable 
 essentially bounded  on $[0,\infty)$ functions;

(a2) $h_k(t)$ are Lebesgue measurable functions,
$$
h_k(t)\leq t, ~\limsup_{t\rightarrow\infty} h_k(t)=\infty;
$$ 

(a3) $\varphi :(-\infty,t_0)\rightarrow R $ is a  Borel
measurable bounded  function.

We assume conditions
(a1)-(a3) hold for all equations throughout the paper.

\noindent
{\bf Definition.} A locally absolutely continuous 
function $x:R\rightarrow R$ is called {\em a solution of the problem}
(\ref{1}), (\ref{2}) if it satisfies the equation (\ref{1})  for almost
all $t\in [t_0,\infty)$ and the equalities (\ref{2}) for $t\leq t_0$.

Below we present a solution representation formula for the nonhomogeneous 
equation with locally Lebesgue integrable right-hand side $f(t)$: 
\begin{equation}
\label{3}
\dot{x}(t)+\sum_{k=1}^m a_k(t)x(h_k(t))=f(t), ~~t \geq t_0.
\end{equation}

\noindent
{\bf Definition}. A solution  $X(t,s)$ of the problem
$$
\dot{x}(t)+ \sum_{k=1}^m a_k(t)x(h_k(t))=0, ~~t\geq s\geq 0,
$$
$$
x(t)=0,~~t< s,~ x(s)=1,
$$
is called {\em the fundamental function} of (\ref{1}).

\newtheorem{uess}{Lemma}
\begin{uess} $\cite{HL,ABR}$
Suppose conditions (a1)-(a3) hold.  
Then the solution of (\ref{3}), (\ref{2})
has the following form
\begin{equation}
\label{4}
x(t)=X(t,t_0)x_0-\int_{t_0}^t X(t,s)\sum_{k=1}^m a_k(s)\varphi(h_k(s))ds
+\int_{t_0}^t X(t,s) f(s)ds,
\end{equation}
where $\varphi(t)=0,~~t\geq t_0$.
\end{uess}

\noindent
{\bf Definition} \cite{HL}. 
Eq.~(\ref{1}) is {\em stable}
if for any initial point $t_0$ and number $\varepsilon>0$ there exists
$\delta>0$ such that the inequality 
$\sup_{t<t_0}|\varphi(t)|+|x(t_0)|<\delta$
implies $|x(t)|<\varepsilon, ~t\geq t_0,$ for the solution 
(\ref{1})-(\ref{2}). 

Eq.~(\ref{1}) is {\em asymptotically stable}
if it stable  
and all solutions of 
(\ref{1})-(\ref{2})  for any initial point $t_0$ 
tend to zero as $t\rightarrow\infty$. 
\vspace{2mm}

In particular, Eq.~(\ref{1}) is asymptotically stable if the fundamental 
function is
uniformly bounded: $|X(t,s)|\leq K, ~t\geq s\geq 0$ 
and all solutions tend to zero as
$t\rightarrow\infty$. 

We apply in this paper only these two conditions of
asymptotic stability.
  
\vspace{2mm}

\noindent
{\bf Definition}. Eq.~(\ref{1}) is 
{\em (uniformly) exponentially stable}, if there exist
$M>0, \mu>0$ such that  the solution of problem (\ref{1})-(\ref{2}) 
has the estimate
$$
|x(t)|\leq M~e^{-\mu (t-t_0)}\left( |x(t_0)|+ 
\sup_{t<t_0}|\varphi(t)| \right),~~t\geq t_ 0,
$$  
where $M$ and $\mu$ do not depend on $t_ 0$.
\vspace{2mm}

\noindent
{\bf Definition}.
The fundamental function  $X(t,s)$ of (\ref{1}) {\em has
an exponential estimation}
if there exist 
$K>0, \lambda>0$ such that 
\begin{equation}
\label{5}
|X(t,s)|\leq K~e^{-\lambda (t-s)},~~t\geq s\geq 0.
\end{equation}

For the linear equation (\ref{1})  with bounded delays the last  two 
definitions are equivalent. For unbounded delays estimation (\ref{5}) 
implies asymptotic stability of
(\ref{1}). 

Under our assumptions the exponential stability does not depend on 
values of equation parameters on any finite interval.

\begin{uess}$\cite{GL,BB4}$
Suppose $a_k(t)\geq 0$. 
If 
\begin{equation}
\label{6}
\int_{\max \{h(t), t_0\}}^t \sum_{i=1}^m a_i(s)ds\leq \frac{1}{e}, ~~ 
h(t)=\min_k \{h_k(t)\}, ~~t\geq t_0,
\end{equation}
or there exists $\lambda >0$, such that
\begin{equation}
\label{7}
\lambda\geq \sum_{k=1}^m A_ke^{\lambda \sigma_k},
\end{equation}
where 
$$
0\leq a_k(t)\leq A_k,~~ t-h_k(t)\leq \sigma_k,~~t\geq t_0,
$$ 
then $X(t,s)>0,~~t\geq s\geq t_0$, where $X(t,s)$ is the fundamental 
function of equation (\ref{1}).
\end{uess}
\begin{uess}$\cite{BB3}$
Suppose $a_k(t)\geq 0$,
\begin{equation}
\label{8}
\liminf_{t\rightarrow\infty}\sum_{k=1}^m a_k(t)>0,
\end{equation}
\begin{equation}
\label{9}
\limsup_{t\rightarrow\infty}(t-h_k(t))<\infty, ~~k=1, \dots , m,
\end{equation}
and there exists $r(t)\leq t$ such that for sufficiently large $t$
$$
\int_{r(t)}^t \sum_{k=1}^m a_k(s)ds\leq\frac{1}{e}.
$$

If
\begin{equation}
\label{10}
\limsup_{t\rightarrow\infty}
\sum_{k=1}^m \frac{ a_k(t)}{\sum_{i=1}^m a_i(t)}
\left|\int_{h_k(t)}^{r(t)} \sum_{i=1}^m a_i(s)ds\right|<1,
\end{equation}
then equation (\ref{1}) is exponentially stable.
\end{uess}
\begin{uess}$\cite{BB3}$
Suppose (\ref{9}) holds and there exists a set of indices $I\subset \{1,\dots,m\}$, such
that  $a_k(t)\geq 0, k\in I$, 
\begin{equation}
\label{10a}
\liminf_{t\rightarrow\infty}\sum_{k\in I}a_k(t)>0,
\end{equation}
and the fundamental function of the equation
\begin{equation}
\label{11}
\dot{x}(t)+\sum_{k\in I}a_k(t)x(h_k(t))=0
\end{equation}
is eventually positive. If
\begin{equation}
\label{12}
\limsup_{t\rightarrow\infty} \frac{\sum_{k\not\in I}
|a_k(t)|}{\sum_{k\in I}a_k(t)}<1,
\end{equation}
then equation (\ref{1}) is exponentially stable.
\end{uess}
The following lemma is a consequence of Corollary 2 \cite{BB5}
obtained for impulsive delay differential equations. 

\begin{uess}
Suppose for equation (\ref{1}) condition (\ref{9}) holds and this equation is exponentially
stable. If
$$
\int_0^{\infty}\sum_{k=1}^n |b_k(s)|ds<\infty, ~\limsup_{t\rightarrow\infty}
(t-g_k(t))<\infty,~ g_k(t)\leq t,
$$ 
then the equation
$$
\dot{x}(t)+\sum_{k=1}^m a_k(t)x(h_k(t))+\sum_{k=1}^n b_k(t)x(g_k(t))=0
$$
is exponentially stable.
\end{uess}

The following elementary result will be used in the paper.

\begin{uess}   
The ordinary differential equation
\begin{equation}
\label{12b}
\dot{x}(t)+a(t)x(t)=0
\end{equation}
is exponentially stable if and only if there exists $R>0$ such that
\begin{equation}
\label{12a}
\liminf_{t\rightarrow\infty}
\int_t^{t+R} a(s)ds >0.
\end{equation}
\end{uess} 

The following example illustrates that a stronger than (\ref{12a}) 
sufficient condition 
\begin{equation}
\label{ode2}
\liminf_{t,s \to \infty} \frac{1}{t-s} \int_s^t
a(\tau)~d\tau >0
\end{equation}
is not necessary for the exponential stability of the ordinary
differential equation (\ref{12b}).
\vspace{2mm}

\noindent
{\bf Example 1.}
Consider the equation
$$
\dot{x}(t)+a(t)x(t)=0, \mbox{~~where~~~}
a(t)= \left\{ \begin{array}{ll} 1, & t \in [2n,2n+1), \\
0, & t \in [2n+1,2n+2),   \end{array}   
\right. ~~n=0,1,2, \dots
$$
Then $\liminf$ in (\ref{ode2}) equals zero, but 
${\displaystyle |X(t,s)|< e~ e^{-0.5(t-s)} }$, so the equation is 
exponentially stable. 
Moreover, if we consider $\liminf$ in (\ref{ode2}) under the 
condition $t-s \geq R$, then it is still zero for any $R \leq 1$.

\section{Main Results}
\label{main}
\begin{uess}
Suppose $a_k(t)\geq 0$, (\ref{8}),  (\ref{9}) hold and
\begin{equation}
\label{13}
\limsup_{t\rightarrow\infty}
\sum_{k=1}^m \frac{ a_k(t)}{\sum_{i=1}^m a_i(t)}
\int_{h_k(t)}^t \sum_{i=1}^m a_i(s)ds<1+\frac{1}{e}.
\end{equation}
Then equation (\ref{1}) is exponentially stable.
\end{uess}   
{\bf Proof.} 
By (\ref{8}) there exists function 
$r(t)\leq t$ such that for sufficiently large $t$
$$
\int_{r(t)}^t \sum_{k=1}^m a_k(s)ds=\frac{1}{e}.  
$$
For this function condition (\ref{10}) has the form
$$
\limsup_{t\rightarrow\infty}
\sum_{k=1}^m \frac{ a_k(t)}{\sum_{i=1}^m a_i(t)}  
\left|\int_{h_k(t)}^t \sum_{i=1}^m a_i(s)ds-
\int_{r(t)}^t \sum_{i=1}^m a_i(s)ds\right|
$$$$
=\limsup_{t\rightarrow\infty}
\sum_{k=1}^m \frac{ a_k(t)}{\sum_{i=1}^m a_i(t)}
\left|\int_{h_k(t)}^t \sum_{i=1}^m a_i(s)ds-\frac{1}{e}\right|<1.
$$
The latter inequality follows from (\ref{13}). 
The reference to Lemma 3 completes the proof. 
\qed

\newtheorem{corollary}{Corollary}
\begin{corollary}
Suppose $a_k(t)\geq 0$, (\ref{8}), (\ref{9}) hold and
\begin{equation}
\label{14}
\limsup_{t\rightarrow\infty}
\int_{\min_k\{h_k(t)\}}^t \sum_{i=1}^m a_i(s)ds<1+\frac{1}{e}.
\end{equation}
Then equation (\ref{1}) is exponentially stable.
\end{corollary}

The following theorem contains stability conditions for equations 
with unbounded delays.
We also omit condition (\ref{8}) in Lemma 7.  
\newtheorem{guess}{Theorem}
\begin{guess}
Suppose $a_k(t)\geq 0$, condition (\ref{13}) holds,
${\displaystyle \sum_{k=1}^m a_k(t)\not= 0}$ a.e. and
\begin{equation}
\label{21}
\int_0^{\infty}\sum_{k=1}^m a_k(t)dt =\infty,~
\limsup_{t\rightarrow\infty} \int_{h_k(t)}^t \sum_{i=1}^m a_i(s)ds <\infty.
\end{equation}
Then equation  (\ref{1}) is asymptotically stable.

If in addition there exists $R>0$ such that
\begin{equation}
\label{22a}
\liminf_{t,\rightarrow\infty}\int_t^{t+R} \sum_{k=1}^ma_k(\tau)d\tau>0
\end{equation}
then the fundamental function of equation  (\ref{1}) has an exponential estimation.

If condition (\ref{9}) also holds then  (\ref{1}) is exponentially stable. 
\end{guess}     
{\bf Proof}.
Let ${\displaystyle s=p(t):=\int_0^t \sum_{k=1}^m a_k(\tau)d\tau,~ y(s)=x(t)}$,
where $p(t)$ is a strictly increasing function.
Then $x(h_k(t))=y(l_k(s)),~ l_k(s)\leq s$, ${\displaystyle l_k(s)=\int_0^{h_k(t)} 
\sum_{k=1}^m a_k(\tau)d\tau}$ 
and (\ref{1}) can be rewritten in the form
\begin{equation}  
\label{22}
\dot{y}(s)+\sum_{k=1}^m b_k(s)y(l_k(s))=0,
\end{equation}
where
${\displaystyle b_k(s)=\frac{a_k(t)}{\sum_{i=1}^m a_i(t)},~~~ 
s-l_k(s)=\int_{h_k(t)}^t \sum_{k=1}^m
a_k(\tau)d\tau.}
$
~~Since ${\displaystyle \sum_{k=1}^m b_k(s)=1}$ and \\
${\displaystyle \limsup_{s\rightarrow\infty}(s-l_k(s))<\infty}$,
then Lemma 7  can be applied to equation (\ref{22}).
We have
$$
\limsup_{s\rightarrow\infty}
\sum_{k=1}^m \frac{ b_k(s)}{\sum_{i=1}^m b_i(s)}
\int_{l_k(s)}^s \sum_{i=1}^m b_i(\tau)d\tau
=\limsup_{s\rightarrow\infty}\sum_{k=1}^m  b_k(s)(s-l_k(s))
$$$$
=\limsup_{t\rightarrow\infty}
\sum_{k=1}^m \frac{ a_k(t)}{\sum_{i=1}^m a_i(t)}
\int_{h_k(t)}^t \sum_{i=1}^m a_i(s)ds<1+\frac{1}{e}.
$$
By Lemma 7  equation (\ref{22}) is exponentially 
stable.
Due to the first equality in (\ref{21})
$t\rightarrow\infty$ implies $s\rightarrow\infty$.
Hence 
${\displaystyle \lim_{t\rightarrow\infty}x(t)=\lim_{s\rightarrow\infty}y(s)=0}$.

Equation (\ref{22}) is exponentially 
stable, thus the fundamental function $Y(u,v)$ of equation (\ref{22}) has 
an exponential estimation
\begin{equation}
\label{22b}
|Y(u,v)|\leq Ke^{-\lambda (u-v)},~u\geq v\geq 0,
\end{equation}
with $K>0,~ \lambda>0$. 
Since ${\displaystyle X(t,s)=Y \left( \int_0^t \sum_{k=1}^m 
a_k(\tau)d\tau,\int_0^s \sum_{k=1}^m a_k(\tau)d\tau \right) }$, where 
$X(t,s)$ is the fundamental function of (\ref{1}), then (\ref{22b}) yields
$$
|X(t,s)|\leq K\exp\left\{ -\lambda \int_s^t \sum_{k=1}^m a_k(\tau)d\tau 
\right\}.
$$
Hence $|X(t,s)|\leq K, ~t\geq s\geq 0$, which together with
${\displaystyle \lim_{t\rightarrow\infty}x(t)=0}$ 
yields that equation (\ref{1}) is asymptotically stable. 

Suppose now that (\ref{22a}) holds. 
Without loss of generality we can assume that for some $R>0, \alpha>0$
we have 
$$\int_t^{t+R} \sum_{k=1}^m a_k(\tau)d\tau\geq \alpha>0,~ t\geq s\geq 0.
$$
Hence
$$
\exp\left\{ -\lambda \int_s^t \sum_{k=1}^m a_k(\tau)d\tau
\right\}\leq \exp\left\{ \lambda R\sup_{t\geq 0}\sum_{k=1}^m
a_k(t)\right\} e^{-\lambda \alpha (t-s)/R}.
$$ 

Thus, condition (\ref{22a}) implies the exponential estimate for $X(t,s)$.

The last statement of the theorem is evident.
\qed

\noindent
{\bf Remark.}
The substitution
${\displaystyle s=p(t):=\int_0^t \sum_{k=1}^m a_k(\tau)d\tau,~ y(s)=x(t)}$
was first used in \cite{LSS}.

\begin{corollary}
Suppose ${\displaystyle a_k(t)\geq 0,~ \sum_{k=1}^m  
a_k(t)\equiv\alpha>0}$, condition (\ref{9}) holds
 and
\begin{equation}
\label{15}
\limsup_{t\rightarrow\infty}
\sum_{k=1}^m  a_k(t)(t-h_k(t))<1+\frac{1}{e}.
\end{equation}
Then equation (\ref{1}) is exponentially stable.  
\end{corollary}
\begin{corollary}
Suppose $a_k(t)=\alpha_k p(t), \alpha_k>0, p(t)>0$ a.e., 
${\displaystyle \int_0^{\infty} p(t)dt =\infty}$ and 
\begin{equation}
\label{16}
\limsup_{t\rightarrow\infty}
\sum_{k=1}^m \alpha_k \int_{h_k(t)}^t p(s)ds<1+\frac{1}{e}.
\end{equation}
Then  equation (\ref{1}) is asymptotically stable.

If in addition there exists $R>0$ such that
\begin{equation}
\label{23a}
\liminf_{t\rightarrow\infty}\int_t^{t+R} p(\tau)d\tau>0,
\end{equation}
then the fundamental function of equation  (\ref{1}) has an exponential estimation.

If also (\ref{9}) holds then equation  (\ref{1}) is exponentially stable.
\end{corollary}
\begin{corollary}  
Suppose $a(t)\geq 0, b(t)\geq 0, a(t)+b(t)\not= 0$ a.e.,
$$
\int_0^{\infty} (a(t)+b(t))dt =\infty,~~~
\limsup_{t\rightarrow\infty}\int_{h(t)}^t (a(s)+b(s))ds<\infty, 
\mbox{~~~and}
$$
\begin{equation} 
\label{19}
\limsup_{t\rightarrow\infty} \frac{b(t)}{a(t)+b(t)}\int_{h(t)}^t
(a(s)+b(s))ds<1+\frac{1}{e}.
\end{equation}
Then the following equation is asymptotically stable 
\begin{equation}
\label{20}
\dot{x}(t)+a(t)x(t)+b(t)x(h(t))=0.
\end{equation}

If in addition there exists $R>0$ such that  
${\displaystyle \liminf_{t\rightarrow\infty}\int_t^{t+R} 
(a(\tau)+b(\tau))d\tau>0}$
then the fundamental function of (\ref{20}) has an exponential estimation.

If also ${\displaystyle \limsup_{t\rightarrow\infty}(t-h(t))<\infty}$
then equation (\ref{20}) is exponentially stable.
\end{corollary}

In the following theorem we will omit condition 
${\displaystyle \sum_{k=1}^m a_k(t)\not= 0}$ a.e.
of Theorem 1. 

\begin{guess}
Suppose $a_k(t)\geq 0$, condition (\ref{14}) and the first inequality 
in (\ref{21}) hold.
Then  equation  (\ref{1}) is asymptotically stable. 

If in addition (\ref{22a}) holds then the fundamental function of 
equation  (\ref{1}) has an exponential estimation.

If also (\ref{9}) holds then equation  (\ref{1}) is exponentially stable.

\end{guess}
{\bf Proof.}
For simplicity suppose that $m=2$ and consider the equation
\begin{equation}
\label{23}
\dot{x}(t)+a(t)x(h(t))+b(t)x(g(t))=0,
\end{equation}
where $a(t)\geq 0, b(t)\geq 0, \int_0^{\infty}(a(s)+b(s))ds=\infty$
and there exist $t_0\geq 0, \varepsilon>0$ such that 
\begin{equation}
\label{24}
\int_{\min\{h(t),g(t)\}}^t (a(s)+b(s))ds<1+\frac{1}{e}-\varepsilon,~t\geq 
t_0.
\end{equation}
Let us find $t_1\geq t_0$ such that $e^{-h(t)}<\varepsilon/4$, 
$e^{-g(t)}<\varepsilon/4$, $t\geq t_1$,
such $t_1$ exists due to (a2). Then $\int_{\min\{h(t),g(t)\}}^t
e^{-s}ds<\varepsilon/2,~t\geq t_1.$ 
Rewrite equation (\ref{23}) in the form
\begin{equation}
\label{25}
\dot{x}(t)+(a(t)+e^{-t})x(h(t))+b(t)x(g(t))-e^{-t}x(h(t))=0,
\end{equation}
where $a(t)+b(t)+e^{-t}>0$. After the substitution 
${\displaystyle s=\int_0^t (a(\tau)+b(\tau)+e^{-\tau})d\tau}$, $y(s)=x(t)$
equation (\ref{25}) has the form
\begin{equation}
\label{26}
\dot{y}(s)+\frac{a(t)+e^{-t}}{a(t)+b(t)+e^{-t}}y(l(s))
+\frac{b(t)}{a(t)+b(t)+e^{-t}}y(p(s))-\frac{e^{-t}}{a(t)+b(t)+e^{-t}}y(l(s))=0,
\end{equation}
where similar to the proof of Theorem 1  
\begin{equation}
\label{27}
s-l(s)=\int_{h(t)}^t (a(\tau)+b(\tau)+e^{-\tau})d\tau,~~~
s-p(s)=\int_{g(t)}^t (a(\tau)+b(\tau)+e^{-\tau})d\tau.
\end{equation}
First we will show that by Corollary 1 the equation
\begin{equation}
\label{28}
\dot{y}(s)+\frac{a(t)+e^{-t}}{a(t)+b(t)+e^{-t}}y(l(s))
+\frac{b(t)}{a(t)+b(t)+e^{-t}}y(p(s))=0
\end{equation}
is exponentially stable. 
Since ~~
${\displaystyle 
\frac{a(t)+e^{-t}}{a(t)+b(t)+e^{-t}}+\frac{b(t)}{a(t)+b(t)+e^{-t}}=1}$,
then (\ref{8}) holds. Condition (\ref{24}) implies (\ref{9}).
So we have to check only condition (\ref{14}) where the sum under the integral is equal to
1.
By (\ref{24}), (\ref{27}) we have 
$$
\int_{\min\{l(s),p(s)\}}^s \!\!\!
1 ds=s-\min\{l(s),p(s)\}, ~~~~~s-l(s)=\int_{h(t)}^t
(a(\tau)+b(\tau)+e^{-\tau})d\tau
$$ $$
=\int_{h(t)}^t (a(\tau)+b(\tau))d\tau
+\int_{h(t)}^t e^{-\tau}d\tau<1+\frac{1}{e}
-\varepsilon+\varepsilon/2=1+\frac{1}{e}-\varepsilon/2, ~t\geq t_1.
$$
The same calculations give $s-p(s)<1+\frac{1}{e}-\varepsilon/2$,
thus condition (\ref{14}) holds. 

Hence equation (\ref{28}) is exponentially stable.

We return now to equation (\ref{26}). We have $ds=(a(t)+b(t)+e^{-t})dt$, 
then
$$
\int_0^{\infty}\frac{e^{-t}}{a(t)+b(t)+e^{-t}}ds=
\int_0^{\infty}\frac{e^{-t}}{a(t)+b(t)+e^{-t}}(a(t)+b(t)+e^{-t})dt<\infty.
$$
By Lemma 5 equation (\ref{26}) is exponentially stable.
Since $t\rightarrow\infty$ implies $s\rightarrow\infty$ then
${\displaystyle 
\lim_{t\rightarrow\infty}x(t)=\lim_{s\rightarrow\infty}y(s)=0}$,
which completes the proof of
the first part of the theorem. The rest of the proof is similar to the 
proof of Theorem 1. 
\qed

\begin{corollary}
Suppose  ${\displaystyle a(t)\geq 0, \int_0^{\infty} a(t)dt =\infty}$ and 
\begin{equation}
\label{17}
\limsup_{t\rightarrow\infty} \int_{h(t)}^t a(s)ds<1+\frac{1}{e}.
\end{equation}
Then the equation 
\begin{equation}
\label{18}
\dot{x}(t)+a(t) x(h(t))=0
\end{equation}
is asymptotically stable.
If in addition condition (\ref{12a}) holds  
then the fundamental function of (\ref{18})
has an exponential estimation.
If also ${\displaystyle \limsup_{t\rightarrow\infty}(t-h(t))<\infty}$
then equation (\ref{18}) is exponentially stable.
\end{corollary}

Now consider equation  (\ref{1}), where 
only some of coefficients are nonnegative.

\begin{guess}
Suppose there exists a set of indices $I\subset \{1,\dots,m\}$ such
that  $a_k(t)\geq 0, k\in I$, 
\begin{equation}
\label{30a}
\int_0^{\infty} \sum_{k\in I}a_k(t)dt =\infty, ~~~
\limsup_{t\rightarrow\infty}\int_{h_k(t)}^t 
\sum_{i\in I}a_i(s)ds<\infty,~~ k=1,\dots,m,
\end{equation}
\begin{equation}
\label{12A}
\sum_{k\not\in I}|a_k(t)|=0, t\in E,~
\limsup_{t\rightarrow\infty, t\not\in E} \frac{\sum_{k\not\in I}
|a_k(t)|}{\sum_{k\in I}a_k(t)}<1, \mbox{~where~}
E=\left\{ t\geq 0, \sum_{k\in I}a_k(t)=0 \right\}.
\end{equation}

If the fundamental function $X_0(t,s)$ of equation  (\ref{11}) is eventually 
positive then all 
solutions of equation  (\ref{1}) tend to zero as
$t\rightarrow\infty$.

If in addition there exists $R>0$ such that 
\begin{equation}
\label{22d}
\liminf_{t\rightarrow\infty}\int_t^{t+R} \sum_{k\in I} a_k(\tau)d\tau>0
\end{equation}
then the fundamental function of equation  (\ref{1}) has an exponential estimation.

If condition (\ref{9}) also holds then  (\ref{1}) is exponentially stable.
\end{guess}
{\bf Proof.}
Without loss of generality we can assume 
$X_0(t,s)>0$, $t\geq s\geq 0$. Rewrite equation (\ref{1}) in the form
\begin{equation}
\label{24a}
\dot{x}(t)+\sum_{k\in I}a_k(t)x(h_k(t))+\sum_{k\not\in I}a_k(t)x(h_k(t))=0.
\end{equation}
Suppose first that ${\displaystyle \sum_{k\in I}a_k(t)\not=0}$ a.e.
After the substitution 
${\displaystyle 
s=p(t):=\int_0^t \sum_{k\in I} a_k(\tau)d\tau}$, $y(s) = x(t)$
we have 
${\displaystyle
x(h_k(t))=y(l_k(s))}$, 
${\displaystyle l_k(s)\leq s,~
l_k(s)=\int_0^{h_k(t)} \sum_{i\in I} a_i(\tau)d\tau, ~k=1,\dots,m,
}$
and (\ref{1}) can be rewritten in the form 
\begin{equation}
\label{24b}
\dot{y}(s)+\sum_{k=1}^m b_k(s)y(l_k(s))=0,
\end{equation}
where
$
{\displaystyle b_k(s)=\frac{a_k(t)}{\sum_{i\in I}  a_i(t)}.}
$
~
Denote by $Y_0(u,v)$ the fundamental function of the equation
$$
\dot{y}(s)+\sum_{k\in I} b_k(s)y(l_k(s))=0.
$$ 
We have
$$X_0(t,s)=Y_0\left( \int_0^t \sum_{k\in I} a_k(\tau)d\tau,\int_0^s
\sum_{k\in I} a_k(\tau)d\tau \right),$$ $$
Y_0(u,v)=X_0(p^{-1}(u),p^{-1}(v))>0, ~u\geq v\geq 0.
$$

Let us check that other conditions of Lemma 4 hold for equation  (\ref{24b}).
Since ${\displaystyle \sum_{k\in I}b_k(s)=1}$ then condition (\ref{10a}) is satisfied.
In addition,
$$
\limsup_{s\rightarrow\infty, p^{-1}(s) \not\in E} \frac{\sum_{k\not\in I}
|b_k(s)|}{\sum_{k\in I}b_k(s)}=
\limsup_{t\rightarrow\infty, t \not\in E} \frac{\sum_{k\not\in I}
|a_k(t)|}{\sum_{k\in I}a_k(t)}<1.
$$
By Lemma 4 equation (\ref{24b}) is exponentially stable. Hence for any solution 
$x(t)$ of (\ref{1}) we have
${\displaystyle \lim_{t\rightarrow\infty}x(t)=\lim_{s\rightarrow\infty}y(s)=0}$.
The end of the proof is similar to the proof of Theorem 2.
In particular, to remove the condition ${\displaystyle \sum_{k\in 
I}a_k(t)\not= 0}$ a.e.
we rewrite the equation by adding the term $e^{-t}$ to one of $a_k(t), 
k\in I$.
\qed

\noindent
{\bf Remark.} Explicit positiveness conditions for the fundamental function
were presented in Lemma 2.

\begin{corollary}
Suppose 
$$
{\displaystyle 
a(t)\geq 0, ~
\int_0^{\infty} a(t)dt =\infty,~ 
\limsup_{t\rightarrow\infty}\int_{g_k(t)}^t a(s)ds<\infty,}
$$$$ 
\sum_{k=1}^n|b_k(t)|=0, ~t\in E,~
{\displaystyle\limsup_{t\rightarrow\infty, t\not\in E}\frac{\sum_{k=1}^n 
|b_k(t)|}{a(t)}<1},
$$
where $E=\{t\geq 0, a(t)=0\}$.
Then the equation 
\begin{equation}
\label{1a}
\dot{x}(t)+a(t)x(t)+\sum_{k=1}^n b_k(t)x(g_k(t))=0
\end{equation}
is asymptotically stable.
If in addition (\ref{12a}) holds 
then the fundamental function of (\ref{1a}) has an exponential estimation.
If also ${\displaystyle \limsup_{t\rightarrow\infty}(t-g_k(t))<\infty}$ 
then (\ref{1a})
is exponentially stable.
\end{corollary}

\begin{guess}
Suppose ${\displaystyle \int_0^{\infty}\sum_{k=1}^m |a_k(s)|ds<\infty}$.
Then all solutions of equation (\ref{1}) are bounded
and (\ref{1}) is not asymptotically stable.
\end{guess}
{\bf Proof.} 
For the fundamental function of (\ref{1}) we have the following estimation
$$
|X(t,s)|\leq \exp\left\{ \int_s^t \sum_{k=1}^m
|a_k(\tau)|d\tau\right\}.
$$
Then by solution representation formula (\ref{4}) for any solution $x(t)$ of (\ref{1}) 
we have
$$
|x(t)|\leq  \exp\left\{ \int_{t_0}^t \sum_{k=1}^m |a_k(s)|ds \right\} |x(t_0)|+
\int_{t_0}^t \exp\left\{ \int_{s}^t \sum_{k=1}^m |a_k(\tau)|d\tau \right\}
\sum_{k=1}^m |a_k(s)||\varphi(h_k(s))|ds
$$$$
\leq \exp\left\{ \int_{t_0}^{\infty} \sum_{k=1}^m |a_k(s)|ds \right\}
\left(|x(t_0)|+\int_{t_0}^{\infty}\sum_{k=1}^m |a_k(s)|ds||\varphi||\right),
$$
where $||\varphi||=\max_{t<0}|\varphi(t)|$.
Then $x(t)$ is a bounded function.

Moreover, ${\displaystyle |X(t,s)| \leq A:= \exp \left\{ \int_{0}^{\infty} 
\sum_{k=1}^m |a_k(s)|ds \right\} }$, $t\geq s \geq 0$. 
Let us choose $t_0 \geq 0$ such that ${\displaystyle \int_{t_0}^{\infty} 
\sum_{k=1}^m |a_k(s)|ds < \frac{1}{2A} }$, then 
${\displaystyle
X'_t(t,t_0)+\sum_{k=1}^m a_k(t) X(h_k(t),t_0)=0,~ X(t_0,t_0)=1
}$ implies 
${\displaystyle
X(t,t_0) \geq 1 - \int_{t_0}^{\infty} \sum_{k=1}^m |a_k(s)|~A~ds
> 1 - A \frac{1}{2A}= \frac{1}{2},}$
thus $X(t,t_0)$ does not tend to zero, so 
(\ref{1}) is not asymptotically stable.
\qed

Theorems 3 and 4  imply the following results.

\begin{corollary}
Suppose $a_k(t)\geq 0$, there exists a set of indices 
$I\subset \{1,\dots,m\}$ such that
condition (\ref{12A}) and the 
second condition in (\ref{30a}) hold. Then all solutions 
of (\ref{1}) are bounded.
\end{corollary}
{\bf Proof.}
If ${\displaystyle \int_0^{\infty}\sum_{k\in I} |a_k(t)|dt=\infty}$,
then all solutions of (\ref{1}) are bounded by Theorem 3.
Let ${\displaystyle \int_0^{\infty}\sum_{k\in I} |a_k(t)|dt<\infty}$.
By (\ref{12}) we have
${\displaystyle \int_0^{\infty}\sum_{k\not\in I} |a_k(t)|dt\leq
\int_0^{\infty}\sum_{k\in I} |a_k(t)|dt<\infty}$.
Then
${\displaystyle \int_0^{\infty}\sum_{k=1}^m |a_k(t)|dt<\infty}$.
By Theorem 4 all solutions of (\ref{1}) are bounded.
\qed

\begin{guess}
Suppose $a_k(t) \geq 0$. If (\ref{1}) is asymptotically 
stable, then the ordinary differential equation 
\begin{equation}
\label{ODE}
\dot{x}(t)+ \left( \sum_{k=1}^ma_k(t) \right) x(t)=0
\end{equation}
is also asymptotically stable.
If in addition (\ref{9}) holds and (\ref{1}) is exponentially stable,
then (\ref{ODE}) is also exponentially stable.
\end{guess}
{\bf Proof.} The solution of (\ref{ODE}), with the initial condition 
$x(t_0)=x_0$, can be presented as \\
${\displaystyle x(t)=x_0 \exp \left\{ -\int_{t_0}^t \sum_{k=1}^ma_k(s)~ds
\right\} }$, so (\ref{ODE}) is asymptotically stable, as far as
\begin{equation}
\label{ode1}
\int_0^{\infty} \sum_{k=1}^ma_k(s)~ds = \infty
\end{equation}
and is exponentially stable if  (\ref{22a}) holds (see Lemma 6).

If (\ref{ode1}) does not hold, then by Theorem 4 equation (\ref{1}) is
not asymptotically stable. 

Further, let us demonstrate that exponential stability of (\ref{1}) 
really implies (\ref{22a}). 

Suppose for the fundamental function of
(\ref{1}) inequality (\ref{5}) holds and condition (\ref{22a})
is not satisfied. 
Then there exists a sequence $\{t_n\}$, 
$t_n\rightarrow\infty$,
such that
\begin{equation}
\label{1star}
\int_{t_n}^{t_n+n} \sum_{k=1}^m a_k(\tau)~d\tau 
< \frac{1}{n}<\frac{1}{e}, ~~~n\geq 3.
\end{equation}
By (\ref{9}) there exists $n_0>3$ such that $t-h_k(t)\leq n_0,~k=1,\dots,m$. 
Lemma 2 implies
$X(t,s)>0$, $t_n\leq s \leq t \leq t_n+n,~n\geq n_0$. 
Similar to the proof of Theorem 4 and using the inequality ${\displaystyle 
1-x \geq e^{-x}}, ~x>0$,
we obtain
$$X(t_n,t_n+n)  \geq 1-\int_{t_n}^{t_n+n} \sum_{k=1}^m a_k(\tau)~d\tau  \geq
\exp\left\{ - \int_{t_n}^{t_n+n}\sum_{k=1}^m a_k(\tau)~d\tau \right\} >
e^{ -\frac{1}{n}}. $$ 
Inequality (\ref{5}) implies $|X(t_n+n,t_n)|\leq Ke^{-\lambda n}$.
Hence $Ke^{-\lambda n}\geq e^{ -\frac{1}{n}},~ n\geq n_0 $,
or $K>e^{\lambda n-1/3}$ for any $n\geq n_0 $.
The contradiction proves the theorem.
\qed

Theorems 3 and 5 imply the following statement.

\begin{corollary}
Suppose $a_k(t)\geq 0$ and the fundamental function of equation (\ref{1})
is positive. Then (\ref{1}) is asymptotically stable if
and only if the ordinary differential equation 
(\ref{ODE}) is asymptotically stable.

If in addition (\ref{9}) holds then (\ref{1}) is exponentially stable
if and only if (\ref{ODE}) is exponentially stable.
\end{corollary}

\section{Discussion and Examples}

In  paper \cite{BB2} we gave a review of known stability tests for 
the linear equation (\ref{1}). In this part we will compare the
new results obtained in this paper with known stability conditions. 

First let us compare the results of the present
paper with our papers \cite{BB1}-\cite{BB3}.
In all these three papers we apply the same method 
based on Bohl-Perron type theorems and
comparison with known exponentially stable equations.

In \cite{BB1}-\cite{BB3} we considered exponential stability only.
Here we also give explicit conditions for asymptotic stability.
For this type of stability, we omit the requirement that the delays are bounded
and the sum of the coefficients is separated from zero.
We also present some new stability tests, based on the results obtained in 
\cite{BB3}.

Compare now the results of the paper with some other known results 
\cite{Y,Kz,SYC,GH2,GD,HL}.
First of all  we replace the constant $\frac{3}{2}$ in most of these
tests by the constant $1+\frac{1}{e}$. Evidently 
${\displaystyle 1+\frac{1}{e}=1.3678\dots<\frac{3}{2}}$,
so we have a worse constant, 
but it is an open problem to obtain $\frac{3}{2}$-stability results for
equations with measurable coefficients and delays.

Consider now equation (\ref{18}) with a single delay. This equation is
well studied beginning with the classical stability result by Myshkis \cite{My}.
We present here 3 statements which cover most of known 
stability tests for this equation.
\vspace{2mm}

\noindent
{\bf Statement 1} $\cite{Y}$.
{\em Suppose $a(t)\geq 0, h(t)\leq t$   
are continuous functions and
\begin{equation}
\label{29a}
\limsup_{t\rightarrow\infty}\int_{h(t)}^t a(s)ds\leq\frac{3}{2}.
\end{equation}
Then all solutions of (\ref{18}) are bounded.

If in addition
$$
\liminf_{t\rightarrow\infty}\int_{h(t)}^t a(s)ds>0,
$$
and the strict inequality in (\ref{29a}) holds
then equation (\ref{18}) is exponentially stable.}
\vspace{2mm}

\noindent
{\bf Statement 2} $\cite{SYC}$.
{\em Suppose $a(t)\geq 0, h(t)\leq t$   
are continuous functions, the strict inequality (\ref{29a})
holds and $\int_0^{\infty} a(s)ds=\infty $. Then all solutions of 
(\ref{18}) tend to zero as $t\rightarrow\infty$.
}
\vspace{2mm}

\noindent
{\bf Statement 3} $\cite{GH2,GD}$.
{\em Suppose $a(t)\geq 0, h(t)\leq t$   
are measurable functions, $\int_0^{\infty}a(s)ds=\infty$,
$A(t)=\int_0^t a(s)ds$ is a strictly monotone increasing function and
$$
\limsup_{t \rightarrow\infty}\int_{h(t)}^t a(s)ds<
\sup_{0<\omega<\pi/2}\left(\omega+\frac{1}{\Phi(\omega)}\right)\approx
1.45\dots,
$$
$\Phi(\omega)=\int_0^{\infty}u(t,\omega)dt$, where $u(t,\omega)$
is a solution of the initial value problem 
$$
y(t)+y(t-\omega)=0,~~ y(t)=0, ~~ t<0, ~~ y(0)=1.
$$
Then equation (\ref{18}) is asymptotically stable.
}
\vspace{2mm}

\noindent
{\bf Example 2.} 
Consider the equation
\begin{equation}
\label{42}
\dot{x}(t)+\alpha(|\sin t|-\sin t)x(h(t))=0, h(t)\leq t, 
\end{equation}
where $h(t)$ is an arbitrary measurable function such that $t-h(t)\leq
\pi$ and $\alpha >0$.

This equation has the form (\ref{18}) where $a(t)=\alpha(|\sin t|-\sin 
t)$.
Let us check that the conditions of  Corollary 5 hold.
It is evident that 
$\int_0^{\infty}a(s)ds =\infty.$
We have
$$
\limsup_{t\rightarrow\infty}\int_{h(t)}^ta(s)ds\leq
\limsup_{t\rightarrow\infty}\int_{t-\pi}^ta(s)ds\leq 
-\alpha\int_{\pi}^{2\pi}2\sin sds=4\alpha.
$$
If $\alpha<0.25\left(1+\frac{1}{e}\right)$,
then condition (\ref{17}) holds, hence all solutions of equation (\ref{42}) 
tend to zero as $t\rightarrow\infty$.

Statements 1-3 fail for this equation.
In Statements 1,2 the delay should be continuous.
In Statement 3 function $A(t)=\int_0^t a(s)ds$ should be strictly increasing.
\vspace{2mm}

Consider now the general equation (\ref{1}) with several delays. 
The following two statements are well known for this equation.
\vspace{2mm}

\noindent
{\bf Statement 4} $ \cite{Kz}$.
{\em Suppose $a_k(t)\geq 0, h_k(t)\leq t$
are continuous functions and
\begin{equation}
\label{42a}
\limsup_{t\rightarrow\infty}a_k(t)
\limsup_{t\rightarrow\infty}(t-h_k(t))\leq 1.
\end{equation}
Then all solutions of (\ref{1}) are bounded 
and 1 in the right hand side of (\ref{42a}) is the best 
possible constant.

If ${\displaystyle \sum_{k=1}^m a_k(t)>0}$ and the strict inequality in (\ref{42a}) is 
valid then all solutions of 
 (\ref{1}) tend to zero as $t\rightarrow\infty$.

If $a_k(t)$ are constants then in (\ref{42a})
the number 1 can be replaced by $3/2$.}
\vspace{2mm}

\noindent
{\bf Statement 5} $ \cite{SYC}$.
{\em Suppose $a_k(t)\geq 0, h_k(t)\leq t$
are continuous, $h_1(t)\leq h_2(t)\leq\dots\leq h_m(t)$
and
\begin{equation}
\label{42b}
 \limsup_{t\rightarrow\infty} \int_{h_1(t)}^t \sum_{k=1}^m a_k(s)ds\leq 3/2.
\end{equation}
Then any solution of (\ref{1}) tends to a constant as $t\rightarrow\infty$.

If in addition ${\displaystyle \int_0^{\infty} \sum_{k=1}^m a_k(s)ds=\infty}$,
then all solutions of (\ref{1}) tend to zero as $t\rightarrow\infty$.          }
\vspace{2mm}

\noindent
{\bf Example 3.}
Consider the equation
\begin{equation}
\label{44}
\dot{x}(t)+\frac{\alpha}{t}x\left( \frac{t}{2}-\sin t 
\right) +\frac{\beta}{t}x\left( \frac{t}{2} \right) =0,~ t\geq
t_0>0, 
\end{equation}
where $\alpha>0, \beta>0$.
Denote 
${\displaystyle
p(t)=\frac{1}{t},~~
h(t)=\frac{t}{2}-\sin t, ~~ g(t)=\frac{t}{2}. }$

We apply Corollary 3. 
Since ${\displaystyle \lim_{t\to \infty}
\left[ \ln\left( \frac{t}{2} \right)- \ln\left( \frac{t}{2}-\sin t\right)
\right]=0}$, then 
$$\limsup_{t\rightarrow\infty}\left(\alpha\int_{h(t)}^t p(s)ds+
\beta\int_{g(t)}^t p(s)ds\right)
\leq (\alpha+\beta)\ln2.$$
Hence if ${\displaystyle \alpha+\beta<\frac{1}{\ln2}\left( 
1+\frac{1}{e} \right)}$ then 
equation (\ref{44}) is exponentially stable.
Statement 4 fails for this equation since the delays are unbounded. 
Statement 5 fails for this equation since neither $h(t)\leq g(t)$ nor $g(t)\leq h(t)$ 
holds.
\vspace{2mm}

Stability results where the nondelay term dominates over the delayed 
terms are well known
beginning with the book of Krasovskii~\cite{Kr}. 
The following result is cited from the monograph~\cite{HL}.
\vspace{2mm}

\noindent
{\bf Statement 6} $\cite{HL}$.
{\em Suppose $a(t), b_k(t), t-h_k(t)$ are bounded continuous functions,
there exist $\delta,k$, $\delta>0, 0<k<1,$ such that
$a(t)\geq \delta$ and  ${\displaystyle \sum_{k=1}^m 
|b_k(t)|<k\delta}$.
Then the equation
\begin{equation}
\label{45a}
\dot{x}(t)+a(t)x(t)+\sum_{k=1}^m b_k(t)x(h_k(t))=0
\end{equation}
is exponentially stable. }
\vspace{2mm}

In Corollary 6 we obtained a similar result without the assumption that the parameters of 
the equation are continuous functions and the delays are bounded.
\vspace{2mm}

\noindent
{\bf Example 4.}
Consider the equation
\begin{equation}
\label{45}
\dot{x}(t)+\frac{1}{t}x(t)+\frac{\alpha}{t}x\left( \frac{t}{2} \right) =0,
~~~ t\geq t_0>0.
\end{equation}
If $\alpha<1$ then by Corollary 6  all solutions of equation (\ref{45})
tend to zero.
The delay is unbounded, thus Statement 6 fails for this equation.
\vspace{2mm}

In \cite{LP} the authors considered a delay autonomous equation with linear and
nonlinear parts, where  the differential equation with the linear part only has 
a positive fundamental function and the linear part dominates over the nonlinear one.
They generalized the early result of Gy\"{o}ri \cite{G} and some results 
of \cite{FH}.

In Theorem 4 we obtained a similar result for a linear nonautonomous equation 
without the assumption that coefficients and delays are continuous. 
\vspace{2mm}

We conclude this paper with some open problems.  
\vspace{2mm}  

\noindent   
{\bf Open Problem 1.}
{\em Prove or disprove that in Corollary 5 the constant 
${\displaystyle 1+\frac{1}{e}}$
can be  changed by the  constant ${\displaystyle \frac{3}{2} }$.}

Note that  all known proofs with the constant $\frac{3}{2}$ apply
methods which are not applicable for equations with measurable 
parameters. 
\vspace{2mm}

\noindent
{\bf Open Problem 2.}
{\em Suppose $a_k(t)\geq 0$, conditions  (\ref{8}), (\ref{9}) hold and equation 
(\ref{1}) is exponentially stable. 

Prove or disprove that for any $b_k(t)$, $0\leq b_k(t)\leq a_k(t)$,
where ${\displaystyle \liminf_{t\rightarrow\infty} \sum_{k=1}^m b_k(t)>0}$, 
the equation 
$$
\dot{x}(t)+\sum_{k=1}^m  b_k(t)x(h_k(t))=0
$$
with the same delays as in (\ref{1}) 
is also exponentially stable.

Obtain similar result for the asymptotic stability.
}
\vspace{2mm}
 
The solution of the following problems would improve Theorems 1 and 5, respectively.
\vspace{2mm}
 
\noindent
{\bf Open Problem 3.}
{\em Suppose (\ref{8}), (\ref{9}) hold and
$$
\limsup_{t\rightarrow\infty}
\sum_{k=1}^m \frac{ |a_k(t)|}{\sum_{i=1}^m a_i(t)}
\int_{h_k(t)}^t \sum_{i=1}^m a_i(s)ds<1+\frac{1}{e}.
$$
Prove or disprove that equation (\ref{1}) is exponentially stable.
}
\vspace{2mm}
 
\noindent
{\bf Open Problem 4.}
{\em Suppose (\ref{1}) is exponentially   stable. Prove or 
disprove that the
ordinary differential equation (\ref{ODE})
is also exponentially (asymptotically)  stable, without restrictions on 
the signs of coefficients $a_k(t) \geq 0$, as in Theorem 5.
}


\begin{thebibliography}{99}
\bibitem{BB1}
L. Berezansky and E. Braverman,
On stability of some linear and nonlinear delay differential equations,
{\em J. Math. Anal. Appl.} {\bf 314} (2006),  391--411.
\bibitem{BB2}
L. Berezansky and E. Braverman,
On exponential stability of  linear differential
equations with several delays,
{\em J. Math. Anal. Appl.} {\bf 324} (2006), 1336--1355.
\bibitem{BB3}
L. Berezansky and E. Braverman,
Explicit stability conditions for  linear differential equations with several delays,
{\em J. Math. Anal. Appl.} {\bf 332} (2007), 246--264.
\bibitem{YS1}
T. Yoneyama and J. Sugie,
On the stability region of scalar delay-differential
equations, {\em J. Math. Anal. Appl.} {\bf 134} (1988), no. 2, 408--425. 
\bibitem{Y}
T. Yoneyama, 
The 3/2 stability theorem for one-dimensional delay-differential equations with unbounded
delay, {\em J. Math. Anal. Appl.} {\bf 165} (1992), 133--143.
\bibitem{Kz}
T. Krisztin, On stability properties for one-dimensional
functional-differential equations,
{\em Funkcial. Ekvac.} {\bf 34} (1991), 241--256.
\bibitem{SYC}
J.\,W.\,H. So, J.\,S. Yu and M.\,P. Chen, Asymptotic stability for scalar
delay differential equations, {\em Funkcial. Ekvac.} {\bf 39} (1996),
1--17.
\bibitem{GH1}
I. Gy\"{o}ri, F. Hartung and J. Turi, Preservation of stability in delay
equations under delay perturbations, {\em J. Math. Anal. Appl.}
{\bf  220} (1998), 290--312.
\bibitem{GH2}
I. Gy\"{o}ri and F. Hartung, Stability in delay perturbed differential and
difference equations, {\em  Topics in
functional differential and difference equations} (Lisbon, 1999),
{\em  Fields Inst. Commun.}, {\bf 29}, Amer. Math. Soc.,
Providence, RI, 2001, 181--194.
\bibitem{GD}
S.\,A. Gusarenko and A.\,I. Domoshnitski\u\i,  Asymptotic and  
oscillation properties of first-order linear scalar
functional-differential equations,
{\em Differential Equations} {\bf 25} (1989), no. 12, 1480--1491.
\bibitem{W}
T. Wang,  Inequalities and stability for a linear scalar 
functional differential equation,  
{\em J. Math. Anal. Appl.} {\bf 298} (2004), no. 1, 33--44.
\bibitem{SY}
J.\,H. Shen and J.\,S. Yu,
 Asymptotic behavior of solutions of neutral differential
equations with positive and negative coefficients, {\em J. Math. Anal. 
Appl.} {\bf 195} (1995), 517--526. 
\bibitem{WL}
X. Wang and L. Liao,
Asymptotic behavior of solutions of neutral differential equations with
positive and negative coefficients, {\em J. Math. Anal. Appl.} 
{\bf 279} (2003), 326--338.
\bibitem{ZW}
Z. Zhang and Z. Wang, 
Asymptotic behavior of solutions of neutral differential
equations with positive and negative coefficients, 
{\em Ann. Differential Equations} {\bf 17} (2001), no. 3, 295--305. 
\bibitem{ZY}
Z. Zhang and J. Yu,  Asymptotic behavior of solutions
of neutral difference equations with positive
and negative coefficients, {\em Math. Sci. Res. Hot-Line} {\bf 2} (1998), 
no. 6, 1--12.
\bibitem{AS}
N.\,V. Azbelev and P.\,M. Simonov, Stability of Differential
Equations with Aftereffect. {\em Stability and Control:
Theory, Methods and Applications}, {\bf 20}. Taylor $\&$ Francis, London, 
2003. 
\bibitem{LTT1}
A. Ivanov, E. Liz and S. Trofimchuk, Halanay inequality, Yorke 3/2
stability criterion, and differential equations with maxima, 
{\em Tohoku Math. J.} (2)
{\bf 54} (2002), 
277--295. 
\bibitem{LTT2}
E. Liz,  V. Tkachenko and S. Trofimchuk,  
A global stability criterion for scalar
functional differential equations,  {\em SIAM J. Math. Anal.} {\bf 35} 
(2003), 596--622. 
\bibitem{T}
X.\,H. Tang, Asymptotic behavior of delay differential
equations with instantaneously terms,  {\em J. Math. Anal. Appl.} {\bf 
302} (2005), no. 2, 342--359. 
\bibitem{Mal1}
V.\,V. Malygina, Some criteria for stability
of equations with retarded argument, 
{\em Differential Equations} {\bf 28} (1992), no. 10, 1398--1405.
\bibitem{Mal2}
V.\,V. Malygina, Stability of solutions of some linear differential
equations with aftereffect, {\em  Russian Math. (Iz. VUZ)} {\bf 37} 
(1993), no. 5, 63--75. 
\bibitem{HL}
J.\,K. Hale and S.\,M.  Verduyn Lunel, Introduction to
Functional Differential equations. {\em Applied Mathematical Sciences}, 
{\bf 99}. Springer-Verlag, New York, 1993.
\bibitem{ABR}
N.\,V. Azbelev, L. Berezansky and  L.\,F. Rahmatullina, 
A linear functional-differential equation of evolution type, 
{\em Differential Equations} {\bf 13} (1977), no. 11, 1331--1339. 
\bibitem{GL}
I. Gy\"{o}ri and G. Ladas,  Oscillation Theory of Delay Differential 
Equations with Applications,
The Clarendon Press, Oxford University Press, New York, 1991. 
\bibitem{BB4}
L. Berezansky and E. Braverman,
On non-oscillation of a scalar delay differential equation, {\em Dynam.
Systems Appl.} {\bf 6} (1997), no. 4, 567--580. 
\bibitem{BB5}
L. Berezansky and E. Braverman,
Preservation of the exponential stability under perturbations of linear delay impulsive
differential equations,
{\em Zeitschrift fur Analysis und ihre Anwendungen} {\bf 14} (1995), 157-175.
\bibitem{LSS}
G. Ladas, Y.G. Sficas, I.P. Stavroulakis,
Asymptotic behavior of solutions of retarded differential equations,
{\em Proc. Amer. Math. Soc.} {\bf 88} (1983), 247--253.
\bibitem{My}
A.\,D. Myshkis, Linear Differential equations with Retarded Argument,
Nauka, Moscow, 1951.
\bibitem{Kr}
N. Krasovskii, 
Stability of Motion, Nauka, Moscow, 1959. Translation, Stanford University Press, 1963.
\bibitem{LP}
E. Liz, M. Pituk,
Exponential stability in a scalar functional differential equation,
{\em Journal of Inequalities and Applications} (2006), Article ID 37195, 
1--10.
\bibitem{G}
I. Gy\"{o}ri,
Interaction between oscillations and global stability in delay differential equations,
{\em Differential and Integral equations} {\bf 3} (1990), 181--200. 
\bibitem{FH}
T. Faria and W. Huang,
Special solutions for linear functional differential equations and asymptotic behaviour,
{\em Differential and Integral equations} {\bf 18} (2005), 337--360.




\end{thebibliography}
\end{document}